\documentclass[letterpaper,oneside,11pt]{article}

\usepackage[USenglish]{babel} 
\usepackage[T1]{fontenc}
\usepackage[ansinew]{inputenc}

\usepackage{lmodern} 

\usepackage{graphicx} 

\usepackage{amsmath}
\usepackage{amsthm}
\usepackage{amsfonts}
\usepackage{mathtools}
\usepackage[hyphens]{url}

\newtheorem{theorem}{Theorem}
\newtheorem*{lemma}{Lemma}

\begin{document}

\pagestyle{empty} 


\title{A characterization of Wilson-Lerch primes}
\author{John Blythe Dobson (j.dobson@uwinnipeg.ca)}

\maketitle

\tableofcontents 
\cleardoublepage 

\pagestyle{plain} 


\begin{abstract}
\noindent
This note presents criteria in terms of Bernoulli numbers for a number to be simultaneously a Wilson prime and a Lerch prime.

\end{abstract}

\section{Notation}

The Fermat quotient $q_p(a) = (a^{p - 1} - 1)/p$. \\
The Wilson quotient $W_p = (\{p-1\}! + 1)/p$. \\
Bernoulli numbers appear in the even index notation of N\"{o}rlund ($B_0=1$, $B_1=-\frac{1}{2}$, $B_2=\frac{1}{6}$, $B_4=-\frac{1}{30}$, $B_3 = B_5 = B_7 = \dots =0$, etc.).

\section{Introduction}

\noindent
Many kinds of special primes can be characterized by the fact that they satisfy a congruence modulo $p^2$ which is only satisfied modulo $p$ by other primes. For example, the Wieferich primes (OEIS A001220) are defined by $2^{p-1} \equiv 1 \pmod{p^2}$, and the Mirimanoff primes (OEIS A014127) by $3^{p-1} \equiv 1 \pmod{p^2}$, these being the best-known examples of the vanishing of the Fermat quotient modulo $p$. Or more generally, special primes may satisfy a congruence modulo $p^n$ which is only satisfied modulo $p^{n-1}$ by other primes; for example, the Wolstenholme primes (OEIS A088164) are defined by $\binom{2p-1}{p-1} \equiv 1 \pmod{p^4}$.

In a similar manner, a Wilson prime (OEIS A007540) is classically defined by the condition that $p$ divides its Wilson quotient; i.\,e. $(p-1)! \equiv -1 \pmod{p^2}$, and a Lerch prime (OEIS A197632) by the condition that $p$ divides its Lerch quotient (see below); i.\,e. $\sum_{a=1}^{p - 1} q_p(a) \equiv W_p \pmod{p^2}$. The Wilson primes $< 2 \cdot 10^{13}$ are 5, 13, and 563 (\cite{CrandallDilcherPomerance}, \cite{Costa}), and the Lerch primes $< 4,496,113$ are 3, 103, 839, and 2237 \cite{Sondow}, with no overlap between the two sequences in the ranges examined. In this note, we present analogous criteria for a prime to possess both of these properties simultaneously. We do not presume that our criteria have anything new to offer from a computational perspective; and considering that the search for Wilson primes has already been carried nearly to the limits of existing means of computation, it is doubtful whether any actual example could be discovered in the foreseeable future. Nonetheless, we present our results in the hope that they may shed some light on the theoretical possibility, or impossibility, of a Wilson-Lerch prime.


\section{The Wilson quotient}
\label{Wilson}

\noindent
The fundamental definition of the Wilson quotient is $W_p = (\{p-1\}! + 1)/p$. In 1899, Glaisher (\cite{Glaisher2}, p.\ 326) proved that (in the modern notation for the Bernoulli numbers)

\begin{equation} \label{eq:Glaisher}
W_p \equiv B_{p - 1} - 1 + \frac{1}{p} \pmod{p},
\end{equation}

\noindent
and generalizations of this will be found in \cite{Lehmer1938}, pp.\ 254--55; \cite{Vandiver1941}, p.\ 578; and \cite{Carlitz1953b}, pp.\ 166--67. Because $W_p$ is an integer for all primes, the right-hand side is $p$-integral, meaning that when written as a reduced fraction, the denominator is not divisible by $p$. Thus the right-hand side has a $p$-adic valuation of at least 0, or equivalently,

\begin{equation} \label{eq:WilsonPrimeCriterionDiminished}
B_{p - 1} - 1 + \frac{1}{p} \equiv 0 \pmod{p^0}.
\end{equation}

\noindent
Clearly the Wilson primes are distinguished by the stricter congruence

\begin{equation} \label{eq:WilsonPrimeCriterion}
B_{p - 1} - 1 + \frac{1}{p} \equiv 0 \pmod{p}.
\end{equation}

\noindent
Now Kummer's congruence for the Bernoulli numbers, as extended by Johnson (\cite{Johnson}, p. 253) to the case where $p-1$ divides the index, gives

\begin{displaymath}
\frac{B_{m(p - 1)} - 1 + \frac{1}{p}}{m(p-1)} \equiv \frac{B_{p - 1} - 1 + \frac{1}{p}}{p-1} \pmod{p} \quad (p > 2),
\end{displaymath}

\noindent
where $m$ may be any positive integer, even a multiple of $p-1$ or of $p$. Since we do not require this theorem in its full generality, for the sake of simplicity we rewrite it with $m=2$:

\begin{equation} \label{eq:KummerJohnson2}
\frac{B_{2p - 2} - 1 + \frac{1}{p}}{2p-2} \equiv \frac{B_{p - 1} - 1 + \frac{1}{p}}{p-1} \pmod{p} \quad (p > 2).
\end{equation}

\noindent
Multiplying throughout by $p-1$, and using (\ref{eq:Glaisher}), we obtain

\begin{displaymath}
\begin{split}
B_{2p - 2} - 1 + \frac{1}{p} & \equiv 2 \left\{ B_{p-1} - 1 + \frac{1}{p} \right\} \\
& \equiv B_{p-1} - 1 + \frac{1}{p} + W_p \pmod{p};
\end{split}
\end{displaymath}

\noindent
in other words,

\begin{equation} \label{eq:Lehmer}
W_p \equiv B_{2p - 2} - B_{p-1} \pmod{p}.
\end{equation}

\noindent
This is a well-known result of Lehmer (\cite{Lehmer1938}, p.\ 355), but we think the derivation from Johnson's supplement to Kummer's congruence is illuminating. For a Wilson prime, clearly we thus have

\begin{equation} \label{eq:WilsonPrimeCriterionJohnson2}
B_{2p - 2} \equiv B_{p-1} \pmod{p}.
\end{equation}

\noindent
While at first glance this may appear to be a pointless reformulation of (\ref{eq:WilsonPrimeCriterion}), the usefulness of this expression will become apparent below. We will also make use of a result of Slavutskii, who rediscovered Johnson's result (\ref{eq:KummerJohnson2}) and extended it to obtain several theorems connecting three Bernoulli numbers with indices divisible by $p-1$, including the following (\cite{Slavutskii}, p.\ 143):

\begin{equation} \label{eq:Slavutskii}
B_{3p - 3} \equiv 3B_{2p - 2} - 3B_{p - 1} + 1 - \frac{1}{p} \pmod{p^2}.
\end{equation}

\noindent
At the risk of belaboring the obvious, we point out that this implies for all primes

\begin{equation} \label{eq:WilsonSlavutskii0}
B_{3p - 3} - 1 + \frac{1}{p} \equiv 0 \pmod{p^0},
\end{equation}

\noindent
and likewise a condition for the Wilson primes equivalent to (\ref{eq:WilsonPrimeCriterionJohnson2}) and analogous in form to (\ref{eq:WilsonPrimeCriterion}),

\begin{equation} \label{eq:WilsonSlavutskii1}
B_{3p - 3} - 1 + \frac{1}{p} \equiv 0 \pmod{p}.
\end{equation}

\noindent
These congruences may be compared with Theorem 2 below.

\section{The Lerch quotient}

\noindent
In 1905, Lerch (\cite{Lerch}, p. 472, eq.\ 4) proved that

\begin{equation} \label{eq:LerchWilson}
\sum_{a=1}^{p - 1} q_p(a) \equiv W_p \pmod{p} \quad (p > 2).
\end{equation}

\noindent
In homage to this important congruence, Jonathan Sondow (\cite{Sondow}, p.\ 3), defined the Lerch quotient,

\begin{displaymath}
\ell_p = \frac{\sum_{a=1}^{p - 1} q_p(a) - W_p }{p},
\end{displaymath}

\noindent
and a Lerch prime as one that divides this quotient; in other words, a prime for which

\begin{equation} \label{eq:LerchPrimeDefinition}
\sum_{a=1}^{p - 1} q_p(a) \equiv W_p \pmod{p^2}.
\end{equation}

\noindent
In a 1953 paper by Carlitz (\cite{Carlitz1953b}, p.\ 166, eq.\ 4.2), (\ref{eq:Glaisher}) and (\ref{eq:LerchWilson}) are combined and partly strengthened to give

\begin{equation} \label{eq:LerchWilsonModP2}
\sum_{a=1}^{p - 1} q_p(a) \equiv B_{p - 1} - 1 + \frac{1}{p} \pmod{p^2} \quad (p > 3).
\end{equation}

\noindent
This supplies an alternate criterion for a Lerch prime, as one satisfying the congruence

\begin{equation} \label{eq:LerchPrimeDefinition2}
W_p \equiv B_{p - 1} - 1 + \frac{1}{p} \pmod{p^2},
\end{equation}

\noindent
which appears in a slightly different notation in Sondow (p.\ 5, eq.\ 6). It will be noted that  (\ref{eq:LerchPrimeDefinition2}) bears the same relation to (\ref{eq:Glaisher}) as (\ref{eq:LerchWilsonModP2}) bears to (\ref{eq:LerchWilson}); i.\,e. each is a $p^2$ variant on a congruence satisfied by all primes. As we shall see below, the closed form (\ref{eq:LerchPrimeDefinition2}) crucially facilitates the comparison of Lerch primes with Wilson primes.

Incidentally, evaluating the left-hand side of (\ref{eq:LerchWilsonModP2}) when the modulus is a higher power of $p$ is a straightforward task, for by the Euler-MacLaurin summation formula, its value is given exactly by

\begin{equation} \label{eq:EulerMacLaurin}
\sum_{a=1}^{p - 1} q_p(a) = - 1 + \frac{1}{p} + \sum_{j=1}^{p} \binom{p}{j} p\thinspace^{j - 2} B_{p - j} \quad (p > 3),
\end{equation}

\noindent
where $B_{p - j}$ vanishes for all even $j$ except $j = p - 1$. This identity, in which the sum in the right-hand side is really just the usual expansion of $\frac{1}{p^2}\left\{ B_p(p) - B_p \right\}$ with the terms reversed, can be used to obtain congruences like (\ref{eq:LerchWilsonModP2}) to any desired precision, though (\ref{eq:LerchWilsonModP2}) is sufficient for our purpose.

\section{Connecting the Wilson quotient with the Lerch quotient}

\noindent
The Wilson quotient is likewise defined by an identity, which --- as pointed out by Lehmer \cite{Lehmer1937} --- traces back to Euler and appears in an independent proof of (\ref{eq:LerchWilson}) given by Beeger (\cite{Beeger1913}, p.\ 83):

\begin{equation} \label{eq:WilsonBeeger}
\begin{split}
W_p & = \frac{1}{p} \cdot \sum_{a=1}^{p-1} (-1)^a \binom{p-1}{a} \left( a^{p-1} - 1 \right) \\
& = \quad \medspace \sum_{a=1}^{p-1} (-1)^a \binom{p-1}{a} \thinspace q_p(a).
\end{split}
\end{equation}

\noindent
The final step of Beeger's proof depends on the well-known result of Lucas (1879) that $\binom{p-1}{a} \equiv  (-1)^a \pmod{p}$ for all $a$ such that $1 \leq a \leq p-1$. The evaluation of the right-hand side of (\ref{eq:WilsonBeeger}) when the modulus is a higher power of $p$ appears to be in general a much more difficult problem. However, as a starting point we can apply the refinement of Lucas's result by Lehmer (\cite{Lehmer1938}, p.\ 360), which states:

\begin{equation} \label{eq:LucasLehmer}
\binom{p-1}{a} \equiv (-1)^a \left\{1 - p H_a + \frac{p^2}{2} H_a^2 - \frac{p^2}{2} H_{a, 2} \right\} \pmod{p^3},
\end{equation}

\noindent
where $H_a$ is the harmonic number $1 + \frac{1}{2} + \frac{1}{3} + \dots + \frac{1}{a}$, and $H_{a, 2}$ is the generalized harmonic number $1 + \frac{1}{2^2} + \frac{1}{3^2} + \dots + \frac{1}{a^2}$. This result, which in its essence can be traced back to Glaisher \cite{Glaisher3}, and which has been extended to the modulus $p^4$ by Z.\,H.\ Sun (\cite{ZHSun2008}, p.\ 285), may be combined with (\ref{eq:WilsonBeeger}) to give the following refinement of Lerch's congruence (\ref{eq:LerchWilson}):

\begin{equation} \label{eq:WilsonGeneralized}
\begin{split}
W_p \equiv \sum_{a=1}^{p-1} q_p(a) - p \cdot \sum_{a=1}^{p-1} H_a \thinspace q_p(a) + \frac{p^2}{2} \cdot \sum_{a=1}^{p-1} H_a^2 \thinspace q_p(a) - \frac{p^2}{2} \cdot \sum_{a=1}^{p-1} H_{a, 2} \thinspace q_p(a) \\
\pmod{p^3}.
\end{split}
\end{equation}

\noindent
So long as $a \le p-1$ it is obvious that $H_a$ and $H_{a, 2}$ are $p$-integral, and so must be the sums containing them. We may thus deduce directly from (\ref{eq:LucasLehmer}) the weaker

\begin{equation} \label{eq:LucasLehmerWeaker}
\binom{p-1}{a} \equiv (-1)^a \left\{1 - p H_a \right\} \pmod{p^2},
\end{equation}

\noindent
and directly from (\ref{eq:WilsonGeneralized}) the weaker

\begin{equation} \label{eq:WilsonGeneralizedWeaker}
W_p \equiv \sum_{a=1}^{p-1} q_p(a) - p \cdot \sum_{a=1}^{p-1} H_a \thinspace q_p(a) \pmod{p^2}.
\end{equation}

\noindent
The sums over products of Harmonic numbers and Fermat quotients in (\ref{eq:WilsonGeneralized}) are relatively intractable, but having recourse to an evaluation of the Wilson quotient by Sun (\cite{ZHSun2000}, pp.\ 210--13) which was obtained by a quite different method, it is known that

\begin{equation} \label{eq:WilsonCongruenceSun}
W_p \equiv \frac{1}{p} - \frac{B_{p-1}}{p-1} + \frac{B_{2p-2}}{2p-2} - \frac{p}{2}\left( \frac{B_{p-1}}{p-1} \right)^2 \pmod{p^2} \quad (p > 3),
\end{equation}

\noindent
where the sum of the first two terms in the right-hand side is congruent to $W_p \pmod{p}$, and the sum of the last two terms is a multiple of $p$. This result, incidentally, establishes that the mod $p^2$ evaluation of the sum in (\ref{eq:LerchWilsonModP2}) in terms of a Bernoulli number has no such simple counterpart in terms of the Wilson quotient. It also allows us to state:

\begin{lemma}
A Lerch prime $p > 3$ is characterized by the congruence
\begin{equation} \label{eq:LerchPrimeDefinition4}
W_p \equiv \frac{B_{2p-2}}{2p} - \frac{B_{p-1}^2}{2p-2} \pmod{p}.
\end{equation}

\begin{proof}
Sun's mod $p^2$ congruence for the Wilson quotient (\ref{eq:WilsonCongruenceSun}) may be combined with the definition of a Lerch prime based on Carlitz's congruence (\ref{eq:LerchPrimeDefinition2}) to give another sufficient condition for a Lerch prime $>3$:

\begin{displaymath} \label{eq:LerchPrime1}
B_{p-1} - 1 + \frac{1}{p} \equiv \frac{1}{p} - \frac{B_{p-1}}{p-1} + \frac{B_{2p-2}}{2p-2} - \frac{p}{2}\left( \frac{B_{p-1}}{p-1} \right)^2 \pmod{p^2},
\end{displaymath}

\noindent
which upon multiplication throughout by $(p-1)/p$ and the cancellation of like terms gives

\begin{displaymath} \label{eq:LerchPrime1b}
B_{p-1} - 1 + \frac{1}{p} \equiv \frac{B_{2p-2}}{2p} - \frac{B_{p-1}^2}{2p-2} \pmod{p}.
\end{displaymath}

\noindent
Glaisher's congruence (\ref{eq:Glaisher}) states that $W_p \equiv B_{p - 1} - 1 + \frac{1}{p} \pmod{p}$, hence the result follows.

\end{proof}

\end{lemma}

\noindent
We can now give:

\begin{theorem}[First Condition for a Wilson prime to be a Lerch Prime]
A prime $p > 3$ is simultaneously a Wilson prime and a Lerch prime if it satisfies the congruence
\begin{equation} \label{eq:WilsonLerchCombinedCriterionOne}
B_{2p - 2} \equiv B_{p-1} \pmod{p^2}.
\end{equation}
\begin{proof}
Setting the left-hand side of (\ref{eq:LerchPrimeDefinition4}) to 0 and multiplying throughout by $2p(p-1)$ gives

\begin{displaymath}
(p-1) \cdot B_{2p-2} \equiv p \cdot B_{p-1}^2 \pmod{p^2}.
\end{displaymath}

\noindent
Substituting the definition of a Wilson prime (\ref{eq:WilsonPrimeCriterion}) in the form $p \cdot B_{p-1} \equiv p - 1 \pmod{p^2}$ into the right-hand side of the above gives

\begin{displaymath}
(p-1) \cdot B_{2p-2} \equiv (p-1) \cdot B_{p-1} \pmod{p^2},
\end{displaymath}

\noindent
and cancelling the common term $p-1$, the result follows.
\end{proof}
\end{theorem}

\noindent
Finally, (\ref{eq:WilsonLerchCombinedCriterionOne}) can be rewritten using only a single Bernoulli number:

\begin{theorem}[Second Condition for a Wilson prime to be a Lerch Prime]
A prime $p > 3$ is simultaneously a Wilson prime and a Lerch prime if it satisfies the congruence
\begin{equation} \label{eq:WilsonLerchCombinedCriterionTwo}
B_{3p - 3} - 1 + \frac{1}{p} \equiv 0 \pmod{p^2}.
\end{equation}
\begin{proof}
Apply the condition (\ref{eq:WilsonLerchCombinedCriterionOne}) to Slavutskii's congruence (\ref{eq:Slavutskii}).

\end{proof}
\end{theorem}

\section{Conclusion}

\noindent
The three congruences (\ref{eq:KummerJohnson2}), (\ref{eq:WilsonPrimeCriterionJohnson2}), and (\ref{eq:WilsonLerchCombinedCriterionOne}), may be seen as forming a progression of increasing stringency, with (\ref{eq:KummerJohnson2}) characterizing primes in general, (\ref{eq:WilsonPrimeCriterionJohnson2}) the Wilson primes, and (\ref{eq:WilsonLerchCombinedCriterionOne}) the Wilson-Lerch primes. The first, though fundamental, has not been traced earlier than Johnson's paper of 1975 \cite{Johnson}, the second has not been traced earlier than Lehmer's paper of 1938 \cite{Lehmer1938}, and the third, at least in regard to the interpretation given to it herein, is believed to be new.

Finally, an example of a progression of congruences where the only notational change is the escalation of the power of $p$ in the modulus may be seen in (\ref{eq:WilsonSlavutskii0}), (\ref{eq:WilsonSlavutskii1}), and (\ref{eq:WilsonLerchCombinedCriterionTwo}), which characterize the primes, the Wilson primes, and the Wilson-Lerch primes, respectively.

\section*{Acknowledgement}

I am grateful to Jonathan Sondow for valuable suggestions as to the improvement of the presentation.

\end{document}